\documentclass{ifacconf}

\usepackage{graphicx}      
\usepackage{natbib}

\usepackage{amsmath} 
\usepackage{amsfonts} 
\usepackage{mathtools} 
\usepackage{subcaption}
\usepackage{multicol}
\usepackage{algorithm}
\usepackage{algpseudocode}

\begin{document}
\begin{frontmatter}

\title{Inexact GMRES Policy Iteration\\ for Large-Scale Markov Decision Processes} 

\thanks[footnoteinfo]{This work was supported by the European Research
Council under the Horizon 2020 Advanced under Grant 787845 (OCAL).}

\author[First]{Matilde Gargiani} 
\author[First]{Dominic Liao-McPherson} 
\author[Second]{Andrea Zanelli}
\author[First]{John Lygeros}

\address[First]{Automatic Control Laboratory (IfA), ETH Zürich, Switzerland (e-mail: \{gmatilde, dliaomc@ethz.ch, jlygeros\}@ethz.ch).}
\address[Second]{Institute for Dynamic Systems and Control (IDSC), ETH, Zurich (e-mail: zanellia@ethz.ch)}

\begin{abstract}                
Policy iteration enjoys a local quadratic rate of contraction, but its iterations are computationally expensive for Markov decision processes (MDPs) with a large number of states. In light of the connection between policy iteration and the semismooth Newton method and taking inspiration from the inexact variants of the latter, we propose \textit{inexact policy iteration}, a new class of methods for large-scale finite MDPs with local contraction guarantees. We then design an instance based on the deployment of GMRES for the approximate policy evaluation step, which we call inexact GMRES policy iteration. Finally, we demonstrate the superior practical performance of inexact GMRES policy iteration on an MDP with 10000 states, where it achieves a $\times 5.8$ and $\times 2.2$ speedup with respect to policy iteration and optimistic policy iteration, respectively.
\end{abstract}
\begin{keyword}
Optimal control; Dynamic programming; GMRES; Inexact semismooth Newton methods.
\end{keyword}
\end{frontmatter}

\section{Introduction}
Stochastic optimal control problems arise in a variety of applications across different fields~\citep{BERTSIMAS19981, optimalcontrolrobotics} and can be compactly expressed in mathematical terms via a recursive functional equation known as the Bellman equation~\citep{Bellman1956}. Dynamic programming (DP) comprises all methods to solve the Bellman equation, such as value iteration (VI), policy iteration (PI) and their variants~\citep{bertsekas2012DPOC}.
Empirical evidence has shown that, among the dynamic programming methods, PI tends to enjoy the fastest rate of convergence. In addition,~\cite{gargiani22} have proved that for finite MDPs PI is an instance of the semismooth Newton method and therefore, by exploiting the structural properties of the Bellman equation, it is possible to conclude local quadratic rate of convergence. Even though PI converges in very few iterations, its time performance degrades rapidly with the size of the state space. In fact, at each iteration PI requires the exact solution of a system of linear equations with dimension equal to the number of states. While the total number of iterations is not dependent on the size of the MDP, the computational complexity of the exact policy evaluation step is strongly dependent on it, diminishing the computational advantages of PI.

An intuitive way to improve the time complexity of PI is to solve the system of linear equations inexactly.
This is the main idea behind optimistic policy iteration (OPI), where the policy evaluation is solved approximately with a finite number of VI steps~\citep{bertsekas2012DPOC}. Variants of this method include, \textit{e.g.}, the deployment of the Gauss-Seidel and mini-batch versions of VI~\citep{gargiani21}. In~\citep{MRKAIC2002517} the author explores the practical performance of variants of OPI where Krylov methods are used for the approximate policy evaluation step instead of VI. The benchmarks show significant performance improvements with respect to PI and OPI for finite MDPs arising from the discretization of stochastic growth models. Variants of OPI are also studied in~\citep{doi:10.1137/110820920, doi:10.1137/100812641} for financial pricing problems. Their numerical examples show that OPI-type methods are generally signiﬁcantly faster in terms of CPU time compared to the full
PI scheme. Finally,~\cite{10.5555/1046920.1088701} study the performance of different OPI-type methods when used in combination with prioritization, partitioning and reordering heuristics.

In light of the connection between policy iteration and the semismooth Newton method and inspired by the inexact variants of the latter, we propose \textit{inexact policy iteration}, a new class of dynamic programming methods (Section~\ref{sec:IPI}). As in OPI, in inexact policy iteration methods the policy evaluation step is carried out only approximately with an iterative solver; however, the number of inner iterations is not fixed \textit{a priori}, but dictated by a stopping condition which depends on the infinity-norm of the Bellman residual function. 
Unlike~\citep{MRKAIC2002517}, we provide a rigorous analysis of the local contraction properties of the methods in this class (Subsection~\ref{sec:TA}). In Subsection~\ref{sec:iGMRES-PI} we design an instance based on the deployment of GMRES~\citep{Saad1986GMRESAG}, which we call inexact GMRES policy iteration (iGMRES-PI), and we also give theoretical and empirical insights on the advantages of GMRES with respect to VI for the approximate solution of the policy evaluation step. Finally, in Section~\ref{sec:NE} we demonstrate the performance superiority of iGMRES-PI on a large-scale MDP with 10000 states versus the only 500 states MDP used in~\citep{MRKAIC2002517}. Section~\ref{sec:problem_setting} is dedicated to the description of the problem setting and the necessary background material.   


%

\section{Problem Setting \& Background}\label{sec:problem_setting}

We consider infinite horizon discounted cost problems for MDPs $\left\{\mathcal{S}, \mathcal{A}, P, g, \gamma   \right\}$ comprising a finite state space $\mathcal{S} = \left\{ 1,\dots, n  \right\}$, a finite action space $\mathcal{A} = \left\{ 1,\dots, m  \right\}$, a transition probability function $P:\mathcal{S}\times \mathcal{A} \times \mathcal{S}\rightarrow [0,1]$ that defines the probability of ending in state $s'$ when applying action $a$ in state $s$, a stage-cost function $g:\mathcal{S}\times\mathcal{A}\rightarrow \mathbb{R}$ that associates to each state-action pair a bounded cost, and a discount factor $\gamma\in(0,1)$. Throughout the paper, with a slight abuse of notation we use $\mathcal{A}(s)$ to denote the nonempty subset of actions that are allowed at state $s$, $p_{ss'}(a) = P(s, a, s')$ for the probability of transitioning to state $s'$ when the system is in state $s$ and action $a\in\mathcal{A}(s)$ is selected with $\sum_{s'\in\mathcal{S}} p_{ss'}(a) = 1$ for all $s\in\mathcal{S}$ and $a\in\mathcal{A}(s)$.

A \textit{deterministic stationary control policy} $\pi: \mathcal{S}\rightarrow \mathcal{A}$ is a function that maps states to actions, with $\pi(s)\in\mathcal{A}(s)$. We use $\Pi$ to denote the set of all deterministic stationary control policies, from now on simply \textit{policies}. At step $t$ of the decision process under the policy $\pi\in\Pi$, the system is in some state $s_t$ and the action $a_t=\pi(s_t)$ is applied. The discounted cost $\gamma^t g(s_t,a_t)$ is accrued and the system transitions to a state $s_{t+1}$ according to the probability distribution $P(s_t, a_t, \cdot)$. This process is repeated leading to the following cumulative discounted cost
\begin{equation}\label{eq: cost of pi}
V^{\pi}(s) = \lim_{T\rightarrow \infty} \mathbb{E}\left[\, \sum_{t=0}^{T-1} \gamma^t g(s_t, \pi(s_t)) \,\,\Big|\,\,s_0=s\right],
\end{equation}
where $\left\{s_0, \pi(s_0), s_1, \pi(s_1),\dots \right\}$ is the state-action sequence generated by the MDP under policy $\pi$ with initial state $s_0$, and the expected value is taken with respect to the corresponding probability measure over the space of sequences. The transition probability distributions induced by policy $\pi$ can be compactly represented by the rows of an $n\times n$ row-stochastic matrix $\left[P^{\pi}\right]_{ss'}=p_{ss'}(\pi(s))$ for all $s,s' \in \mathcal{S}$ and the costs induced by policy $\pi$ by the vector $g^{\pi} = \begin{bmatrix}
g(1, \pi(1)), &\dots &,  g(n, \pi(n)) 
\end{bmatrix}^\top\in\mathbb{R}^{n}$.
The optimal cost is defined as 
\begin{equation}\label{eq: optimal cost}
V^*(s) = \min_{\pi\in\Pi} V^{\pi}(s)\,,\quad\forall s \in \mathcal{S}. 
\end{equation}
Any policy $\pi^*\in\Pi$ that attains the optimal cost is called an optimal policy. 
Notice that in~\eqref{eq: optimal cost} we restrict our attention to stationary deterministic policies as in our setting there exists a policy in this class that attains $V^*$~\citep[Section 1.1.4]{bertsekas2012DPOC}.

Equations~\eqref{eq: cost of pi} and~\eqref{eq: optimal cost} admit \textit{recursive} formulations which are known as the \textit{Bellman equations}. In particular
\begin{equation}\label{eq: Bellman_eq_pi}
V^{\pi}(s) = g(s, \pi(s)) + \gamma \sum_{s'\in \mathcal{S}}   p_{ss'}(a) V^{\pi}(s')\quad \forall s\in\mathcal{S}\,,
\end{equation}
is the Bellman equation associated with policy $\pi$, and 
\begin{equation}\label{eq: Bellman_eq_star}
V^*(s) = \min_{\pi\in\Pi}\left\{ g(s, \pi(s)) + \gamma \sum_{s'\in \mathcal{S}}   p_{ss'}(a) V^{*}(s') \right\} \,\, \forall s\in\mathcal{S}\,,
\end{equation}
is the Bellman equation associated with the optimal cost.

Given the cost $V:\mathcal{S}\rightarrow \mathbb{R}^n$, any policy which satisfies the following equation
\begin{equation}
\pi(s)\in \arg\min \left\{ g(s, \pi(s)) + \gamma  \sum_{s'\in \mathcal{S}}  p_{ss'}(a) V(s')\right\} \,\,\forall s\in\mathcal{S}
\end{equation}
is called \textit{greedy} with respect to the cost $V$. We denote with $\text{GreedyPolicy}(V)$ the operator which extracts a greedy policy associated with $V$.

Starting from the Bellman equations we can define two mappings, $T^{\pi}: \mathbb{R}^n \rightarrow \mathbb{R}^n$ and $T: \mathbb{R}^n \rightarrow \mathbb{R}^n$, where $T^{\pi}V = g^{\pi} + \gamma P^{\pi}V$ and $TV = \min_{\pi\in\Pi}\left\{ g^{\pi} + \gamma P^{\pi} V \right\}$.
These mappings are known as the \textit{Bellman operators} and allow one to rewrite~\eqref{eq: Bellman_eq_pi} and~\eqref{eq: Bellman_eq_star} in a compact form, $V^{\pi} = T^{\pi}V^{\pi}$ and $V^* = TV^*$, respectively. The Bellman operators are $\gamma$-contractive, monotone and shift-invariant and, in our setting, have $V^{\pi}$ and $V^*$ as their unique fixed-points, respectively. We refer to~\citep{bertsekas2012DPOC} for a detailed discussion on the properties of the Bellman operators.

\subsection{Dynamic Programming}

DP comprises the methods to solve~\eqref{eq: Bellman_eq_star}~\citep{bertsekas2012DPOC}. In this work we focus on variants of PI and VI. In particular, given an arbitrary initial cost vector $V_0\in \mathbb{R}^n$, VI is the fixed-point iteration\newline
\begin{minipage}{0.5\linewidth}
    \begin{equation}
      V_{k+1} = T^{\pi}V_k\,, \label{eq: VI_policyevaluation}
    \end{equation}
\end{minipage}%
\begin{minipage}{0.5\linewidth}
    \begin{equation}
    V_{k+1} = TV_k\,, \label{eq: VI_optimalcost}
    \end{equation}
\end{minipage}

and enjoys global linear convergence to $V^{\pi}$ and $V^*$, respectively, with a $\gamma$-contraction rate.
PI instead starts with an arbitrary policy $\pi_0\in\Pi$ and alternates two steps: \textit{policy evaluation}
\begin{equation}
V^{\pi_k} = \left(I - \gamma P^{\pi_k} \right)^{-1}g^{\pi_k}\,,
\end{equation}
and \textit{policy improvement}
\begin{equation}
\pi_{k+1} \in \arg\min_{\pi\in\Pi} \left\{g^{\pi} + \gamma P^{\pi}V^{\pi_k}\right\}\,.
\end{equation}
Exact PI converges in a finite number of iterations, but the worst-case upper bound for large state spaces could be dramatic. Fortunately, PI enjoys global linear convergence to $V^*$ with rate $\gamma$. 
In addition, the convergence rate superiority of PI with respect to VI has been long suggested by extensive empirical evidence and only recently proved for this setting. In particular,~\cite{gargiani22} show that the solution of the Bellman equation~\eqref{eq: Bellman_eq_star} can be expressed as the root of the so-called \textit{Bellman residual function} $r:\mathbb{R}^n \rightarrow \mathbb{R}^n$ with
\begin{equation}\label{eq:Bellman_res}
r(V) = V - TV\,.
\end{equation}
Consequently, solving the Bellman equation corresponds to computing the root of the Bellman residual function. Applying a semismooth variant of Newton's method to~\eqref{eq:Bellman_res} yields the iteration
\begin{equation}\label{eq:Newton_PIiterate}
V_{k+1} = V_k - J_k^{-1}r(V_k)\,,
\end{equation}
where $J_k$ is an element in Clarke's generalized Jacobian of $r$ at $V_k$.
The authors of~\citep{gargiani22} show that the PI iterate is an instance of~\eqref{eq:Newton_PIiterate} and that the policy evaluation step corresponds to the solution of the Newtonian linear system
\begin{equation}\label{eq:Newtonian_linearsystem}
r(V_k) + J_k \left( V - V_k \right) = 0\,,
\end{equation}
where $J_k = I - \gamma P^{\pi_k}$. This equivalence, together with the structural properties of the Bellman residual function, proves local quadratic convergence of PI. 

Despite its fast convergence rate, PI iterations are computationally expensive for MDPs with a large number of states. In particular, in the scenarios where $n$ is large the policy evaluation step is not practical as it requires the exact solution of an $n$-dimensional system of linear equations. An alternative is OPI, where the policy evaluation step is carried out approximately with a fixed number $W$ of VI steps. This number is generally selected to trade-off computational complexity and convergence rate. Notice that when $W=1$ we resort to VI and when $W\rightarrow \infty$ we resort to PI. We refer to~\citep[Chapter 2]{bertsekas2012DPOC} for a thorough analysis of VI, PI and OPI.

Even if the semismooth Newton method enjoys a fast rate of convergence, computing the exact solution of~\eqref{eq:Newtonian_linearsystem} using a direct method can be expensive if the number of unknowns is large. A more computationally efficient solution in the large-scale case consists in solving~\eqref{eq:Newtonian_linearsystem} only approximately with some iterative linear solver and using a certain stopping rule. These are the principles behind \textit{inexact semismooth Newton methods}~\citep{izmasolo14, MARTINEZ1995127}. In particular, $V_{k+1}$ is no longer required to exactly solve~\eqref{eq:Newtonian_linearsystem}, but only to satisfy 
\begin{equation}\label{eq: stopping_rule}
\Vert r(V_k) + J_k \left( V_{k+1} - V_k \right) \Vert \leq \alpha_k\Vert r(V_k) \Vert
\end{equation}
for some $\alpha_k \in [0,1)$. The sequence $\left\{ \alpha_k\right\}$ is called \textit{forcing sequence} and it greatly affects both local convergence properties and robustness of the method~\citep{izmasolo14}.
Different iterative linear solvers can be used to approximately solve~\eqref{eq:Newtonian_linearsystem}~\citep{hackbusch1994iterative}. Often Krylov subspace methods, such as the generalized minimal residual method (GMRES)~\citep{Saad1986GMRESAG}, are deployed in large-scale scenarios.



\subsection{GMRES}\label{sec: GMRES}

Consider a general system of linear equations 
\begin{equation}\label{eq:linear_system}
A x = b\,,
\end{equation}
where $b\in\mathbb{R}^n$ and $A\in\mathbb{R}^{n\times n}$ is a non-singular matrix. Starting from an initial guess $x_0  \in \mathbb{R}^n$ with residual $\Phi(x_0) = b - Ax_0$, GMRES~\citep{Saad1986GMRESAG} generates a sequence $\left\{ x_i \right\}$ of approximate solutions to~\eqref{eq:linear_system} with 
\begin{equation}\label{eq:GMRES_LSQ}
x_i = \arg\min_{x} \left\{ \Vert b - Ax\Vert_2 : x \in x_0 + \mathcal{K}_i \right\}\,,
\end{equation} 
where $\mathcal{K}_i = \text{span}\left\{ \Phi_0, A\Phi_0, A^2 \Phi_0, \dots, A^{i-1}\Phi_0 \right\}$ is known as the $i$-th Krylov subspace and $\Phi_0 = \Phi(x_0)$. 
In particular, at each iteration GMRES generates an orthonormal basis of $\mathcal{K}_i$ via the \textit{Arnoldi's method}~\citep{Saad1986GMRESAG} and then deploys it to solve~\eqref{eq:GMRES_LSQ}.
Unlike the conjugate gradient method, the orthonormal bases can not be computed with a short recurrence. When $i$ increases the number of stored vectors increases like $i$ and the number of multiplication like $0.5 i^2 n$. A practical variant of GMRES, denoted as GMRES$(i)$, consists in restarting the algorithm after every $i$ iterations.

\begin{figure}
\begin{center}
\includegraphics[width=0.38\textwidth]{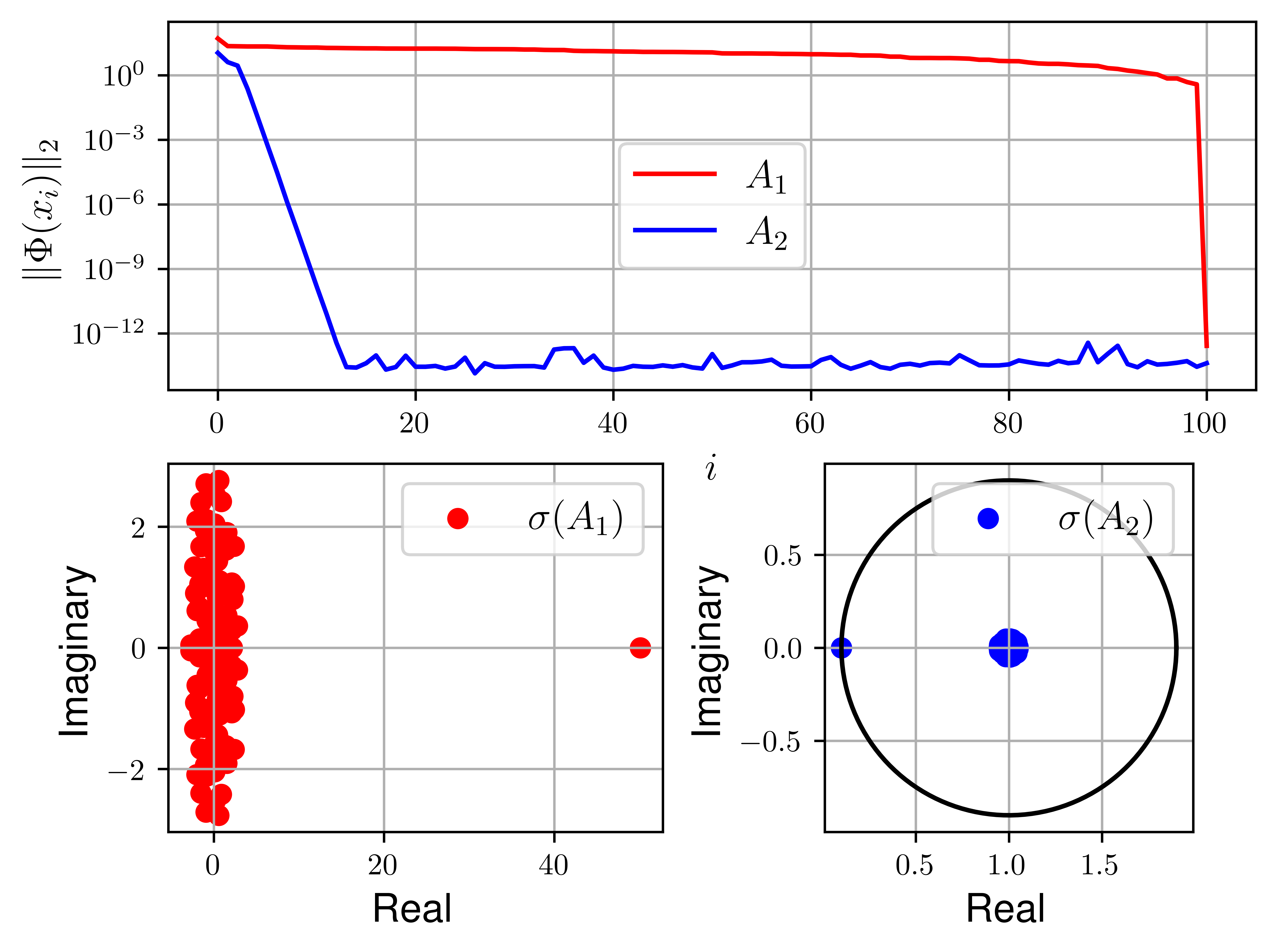}    
\vspace*{-3mm}
\caption{Convergence of GMRES for the case of non-clustered (red) and clustered (blue) eigenvalues of the coefficient matrix.} 
\label{fig:GMRES_convergence}
\end{center}
\end{figure}

GMRES with exact arithmetic converges to the solution of~\eqref{eq:linear_system} in at most $n$ steps. Its convergence rate though is greatly affected by the distribution of the eigenvalues of the coefficient matrix~\citep{GMRES_convclusters}. 
This is exemplified in Figure~\ref{fig:GMRES_convergence}, where GMRES is used to solve the linear systems $A_1 x= b$ and $A_2 x = b$. In particular, $A_1\in\mathbb{R}^{100\times 100}$ is a matrix with non-clustered complex eigenvalues, while all the eigenvalues of $A_2 \in \mathbb{R}^{100\times 100}$ are contained in the circle of center $(1,0)$ and radius $0.9$ in the complex plane. In the first scenario the norm of the residual is significantly decreased only when $i=100$, while in the second scenario we observe R-linear convergence with a fast rate starting from the first iteration.

We refer to~\citep{Saad1986GMRESAG, GMRES_convclusters} for a detailed description of GMRES and its convergence properties. See~\citep[Algorithm 3]{Saad1986GMRESAG} for a pseudocode description of GMRES.

\section{Inexact Policy Iteration Methods}\label{sec:IPI}

We define a novel variant of PI for large-scale scenarios, which we call inexact policy iteration. This class of methods is based on approximately solving the policy evaluation step with an iterative linear solver.
The methods in this class start with an initial guess of the optimal cost $V_0\in\mathbb{R}^n$  and then at every iteration extract a greedy policy associated with the current iterate $V_k\in\mathbb{R}^n$, which is used to compute an element in Clarke's generalized Jacobian.
The next iterate $V_{k+1}\in\mathbb{R}^n$ is selected as an approximate solution of the Newtonian linear system
\begin{equation}\label{eq:IPI_newtonian_ls}
\left( I -\gamma P^{\pi_k} \right) V = g^{\pi_k} 
\end{equation}
which verifies the stopping condition in~\eqref{eq: stopping_rule} with the infinity-norm. 
Because of the specific structure of the Bellman residual function,~\eqref{eq: stopping_rule} simplifies to
\begin{equation*}
\begin{aligned}
\Vert g^{\pi_k} -  \left( I - \gamma P^{\pi_k}\right) V_{k+1} \Vert \leq \alpha_k \Vert g^{\pi_k} -  \left( I - \gamma P^{\pi_k}\right) V_{k}  \Vert\,.
\end{aligned}
\end{equation*}
In principle, any iterative solver for linear systems with non-singular coefficient matrices can be used to generate an approximate solution of~\eqref{eq:IPI_newtonian_ls}, such as VI, its mini-batch version~\citep{gargiani21} and GMRES. 
Notice that, when VI is deployed as inner solver, we obtain a variant of OPI where the number of inner iterations is not selected a priori, but dictated by the stopping condition. 
See Algorithm~\ref{alg:inexact_pi} for a pseudocode description of a general inexact policy iteration method.
\begin{algorithm}
  \caption{Inexact Policy Iteration}\label{alg:inexact_pi}

  \begin{algorithmic}[1]
  \State \textbf{Initialization:}   $V_0 \in \mathbb{R}^n$, $\alpha\in (0,1)$, $K > 0 $

  \For{$k=0,1,\dots, K-1$}
  \State $\pi_k \leftarrow \text{GreedyPolicy}(V_k)$
  \State $J^{\pi_k} = \left( I - \gamma P^{\pi_k} \right)$
  \State $V_{k+1} \leftarrow V_k$
  \While{$\Vert g^{\pi_k} - J^{\pi_k} V_{k+1}  \Vert_{\infty} > \alpha\Vert g^{\pi_k} - J^ {\pi_k} V_k \Vert_{\infty}$}
   
   \State $V_{k+1} \leftarrow \text{IterativeLinearSolver}\left( J^{\pi_k}, g^{\pi_k}, V_{k+1}\right)$
  
  \EndWhile  	
  \EndFor
  
  \end{algorithmic}
\end{algorithm}

\subsection{Inexact GMRES Policy Iteration}\label{sec:iGMRES-PI}
We deploy the presented algorithmic framework to design a novel DP method for large-scale applications. 
The selection of the inner solver is important for the performance of the overall scheme, as a more efficient solver will require less time to meet the stopping condition, leading to an overall faster method.
\begin{figure}
\begin{center}
\includegraphics[width=0.38\textwidth]{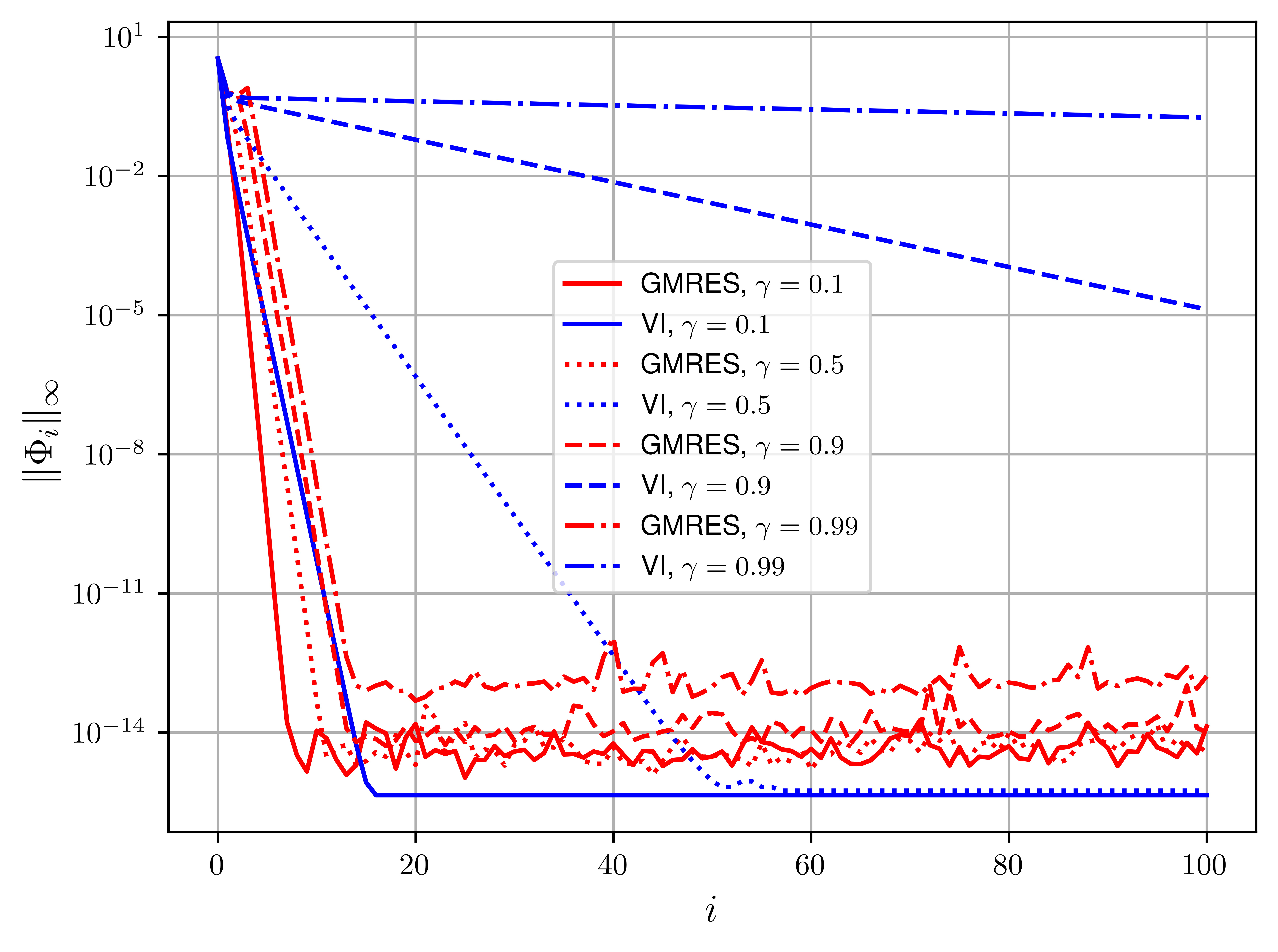}    
\vspace*{-3mm}
\caption{Convergence of GMRES (red) and VI (blue) for policy evaluation for different values of $\gamma$ and $n=100$.} 
\label{fig:GMRESvsVI_convergence}
\end{center}
\end{figure}
Given the particular structure of the Newtonian linear system in~\eqref{eq:IPI_newtonian_ls}, we propose to deploy GMRES as iterative linear solver in Step 7 of Algorithm~\ref{alg:inexact_pi}. In particular, the coefficient matrices are non-singular~\citep[Proposition 3.3]{gargiani22}, but, unless stronger assumptions on the geometry of the underlying MDP hold, we can not rely on symmetry. In addition, as discussed in Section~\ref{sec: GMRES}, GMRES has a particularly favorable convergence behavior for the case of coefficient matrices with clustered eigenvalues. The following lemma demonstrates that this is exactly the scenario encountered in inexact policy iteration methods as the eigenvalues of the coefficient matrices of the Newtonian linear systems are clustered in a circle of radius less than 1.

\begin{lem}
For any $\pi \in \Pi$, the eigenvalues of $I-\gamma P^{\pi}$ are contained in a circle centered at $(1, 0)$ and with radius $\gamma$ in the complex plane.
\end{lem}
\begin{pf}
The result follows directly from the fact that the eigenvalues of $P^{\pi}$ are contained in a circle centered at $(0,0)$ and with radius $1$ in the complex plane~\citep[Theorem 5.3 in Chapter 2]{doi:10.1137/1.9781611971262}.
\end{pf}
Finally, as depicted in Figure~\ref{fig:GMRESvsVI_convergence} for policy evaluation, the convergence rate of GMRES appears to be more robust against the discount factor compared to that of VI.
Algorithm~\ref{alg:gmres_pi_restarted} provides a pseudocode description of iGMRES-PI, where we have adopted the restarted version of GMRES to reduce the computational and storage complexity.




\begin{algorithm}
  \caption{$\text{GMRES}(W)$ Policy Iteration}\label{alg:gmres_pi_restarted}

  \begin{algorithmic}[1]
  \State \textbf{Initialization:} $V_0 \in \mathbb{R}^n$, $\alpha \in (0, 1)$, $K, W > 0$

  \For{$k = 0,\dots, K-1$}
       \State $\pi_k \leftarrow \text{GreedyPolicy}(V_{k})$
	   \State $J^{\pi_k} \leftarrow \left(I-\gamma P^{\pi_k} \right)$  
	   \State $V_{k+1} \leftarrow V_{k}$     
	   \State $\Phi_k\leftarrow g^{\pi_k} - J^{\pi_k}V_k$       
       \State $Q \leftarrow \mathbf{0}_{n\times (W+1)}$
       \State $H \leftarrow \mathbf{0}_{(W+1)\times W}$
       \State $[Q]_{\cdot 1} \leftarrow \Phi_k/\Vert \Phi_k \Vert_2$
	   \State $i\leftarrow 1$       
       \While{$\Vert g^{\pi_k} - J^{\pi_k} V_{k+1} \Vert_{\infty} > \alpha \Vert \Phi_k \Vert_{\infty}$}
       \State $q \leftarrow J^{\pi_k} [Q]_{\cdot i}$
       \For{$j=1,\dots, i$}
       \State $[H]_{ji} \leftarrow [Q]_{\cdot j}^{\top}q$
       \State $q \leftarrow q - [H]_{ji}[Q]_{\cdot j}$
       \EndFor
       \State $[H]_{i+1, i} \leftarrow \Vert q \Vert_2$
       \If{$\Vert q \Vert_2=0$}
       \State $\tilde{y} \leftarrow \arg\min_{y} \Vert \Vert \Phi \Vert_2 \cdot e_1 - Hy \Vert_2$
       \State $V_{k+1} \leftarrow Q \tilde{y} + V_k$
       \State{go to line 35}
       \EndIf		  
       \State $[Q]_{\cdot, i+1} \leftarrow q/ [H]_{i+1, i}$
       \State $\tilde{y} \leftarrow \arg\min_{y} \Vert \Vert \Phi \Vert_2 \cdot e_1 - Hy \Vert_2$
       \State $V_{k+1} \leftarrow Q \tilde{y} + V_k$
       \If{$i=W$}
       \State $Q \leftarrow \mathbf{0}_{n\times (W+1)}$
       \State $H \leftarrow \mathbf{0}_{(W+1)\times W}$
       \State $[Q]_{\cdot 1} \leftarrow \left(g^{\pi_k} - J^{\pi_k}V_{k+1}\right)/\Vert g^{\pi_k} - J^{\pi_k}V_{k+1} \Vert_2$
       \State $i\leftarrow 1$
       \Else 
       \State $i\leftarrow i+1$	
       \EndIf
      
       \EndWhile 
  \EndFor
  \end{algorithmic}
\end{algorithm}

\subsection{Theoretical Analysis}\label{sec:TA}

In this section we provide an analysis of the local convergence properties of inexact policy iteration methods for finite MDPs with discount factor $\gamma \in (0,1)$. We start by characterizing the Lipschitz constant of the Bellman residual function and deriving an upper bound on the infinity-norm of the inverse of the coefficient matrix of the Newtonian linear system in~\eqref{eq:IPI_newtonian_ls}. 

\begin{lem}\label{lm:lipschitz_const}
Let $r:\mathbb{R}^n \rightarrow \mathbb{R}^n$ be the Bellman residual function as defined in~\eqref{eq:Bellman_res}. Then,
\begin{equation*}
\Vert r(V_1) - r(V_2) \Vert_{\infty} \leq (1+\gamma)\Vert V_1 - V_2 \Vert_{\infty}\,,\quad \forall \,\, V_1, V_2 \in \mathbb{R}^n\,.
\end{equation*} 
\end{lem}
\begin{pf}
As shown in~\citep{gargiani22}, $r$ is piecewise affine with selection functions $r_{\pi}(V) = V - T^{\pi}V = (I-\gamma P^{\pi})V - g^{\pi}$ for all $\pi\in\Pi$. Piecewise affine functions are globally Lipschitz continuous and their Lipschitz constant is given by the maximum over the norms of the coefficient matrices of their selection functions~\citep[Proposition 4.2.2]{facchineipangvol12003}. 
Therefore, by exploiting the fact that for any $\pi\in\Pi$ the matrices $P^{\pi}$ are row-stochastic, we obtain 
\begingroup
\allowdisplaybreaks
\begin{align*}
&\max_{\pi\in\Pi}\Vert I - \gamma P^{\pi} \Vert_{\infty} \\
&=\max_{\pi\in\Pi} \max_{s\in\mathcal{S}} \!\!\sum_{\substack{s'\in\mathcal{S}\setminus \left\{ s\right\}}}\!\!\!\! \vert -\gamma p_{ss'}(\pi(s))\vert + \vert 1 - \gamma p_{ss}(\pi(s))\vert\\
&= \max_{\pi\in\Pi} \max_{s\in\mathcal{S}} \,\gamma \left(1 - p_{ss}(\pi(s)) \right) + 1 - \gamma p_{ss}(\pi(s))\\
&\leq 1+\gamma\,,
\end{align*}
\endgroup
which concludes the proof.~\qed
\end{pf}
\begin{lem}\label{lm:inverse_bound}
For any $\pi\in \Pi$ the following inequality holds
\begin{equation*}
\Vert \left(I - \gamma P^{\pi}\right)^{-1} \Vert_{\infty} \leq \frac{1}{1-\gamma}\,.
\end{equation*}
\end{lem}
\vspace*{-4mm}
\begin{pf}
Since $\Vert \gamma P^{\pi} \Vert_{\infty} = \gamma <1$ for any $\pi\in\Pi$, then $I - \gamma P^{\pi}$ is invertible and $\left( I - \gamma P^{\pi}\right)^{-1} = \sum_{k=0}^{\infty} \left( \gamma P^{\pi}\right)^k$~\citep[Chapter 10]{Axler2020}. Therefore for any $\pi\in \Pi$
\begin{equation*}
\begin{aligned}
\Vert \left( I - \gamma P^{\pi} \right)^{-1} \!\Vert_{\infty}&=  \Big\Vert \sum_{k=0}^{\infty}\left( \gamma P^{\pi}\right)^k \Big\Vert_{\infty}\\
&\leq \sum_{k=0}^{\infty} \Vert \gamma P^{\pi} \Vert_{\infty}^k\\ 
&= \frac{1}{1-\gamma}\,,  
\end{aligned}
\end{equation*}
where the first inequality follows from the properties of the infinity-norm and the last equality follows from the properties of the geometric series.~\qed
\end{pf}

The following theorem characterizes the local contraction of inexact policy iteration methods. 
\begin{thm}[local contraction]\label{th:local_contraction}
Consider a general inexact policy iteration method as given in Algorithm~\ref{alg:inexact_pi}. Assume that $\pi_k$ in Step 3 is a non-spurious greedy policy~\citep[Definition~3.2]{gargiani22} and let $\left\{ \alpha_k \right\}$ be a sequence of positive numbers contained in $[0, \alpha]$, with $\alpha \in (0, \frac{1-\gamma}{1+\gamma})$. Then there exists a neighborhood of $V^*$ such that, for any $V_0\in \mathbb{R}^n$ in this neighborhood, the inexact policy iteration method is Q-linearly convergent to $V^*$ with rate $\frac{1+\gamma}{1-\gamma}\alpha$. If $\lim_{k\rightarrow \infty} \alpha_k = 0$, then the method enjoys local Q-superlinear convergence. 
\end{thm}
\begin{pf}
Since by assumption $\pi_k$ is a non-spurious greedy policy, then $I-\gamma P^{\pi_k}$ is an element in Clarke's generalized Jacobian of $r$ at $V_k$~\citep{gargiani22}. In addition, since $r$ is globally CD-regular~\citep[Proposition~3.3]{gargiani22}, the sequence~\eqref{eq:Newton_PIiterate} is globally well-defined.
In the following derivations we use $J_{k} = I-\gamma P^{\pi_k}$ and $\Delta V_{k} = V_{k+1} - V_k$ 
\begingroup
\allowdisplaybreaks
\begin{align*}
&\Vert V_{k+1} - V^* \Vert_{\infty} = \Vert V_k + \Delta V_{k} - V^* \Vert_{\infty}\\
&= \Vert V_k + J_{k}^{-1}J_{k} \Delta V_{k} - V^* \Vert_{\infty}\\
&= \Vert V_k - J_{k}^{-1}r(V_k) + J_{k}^{-1}r(V_k) + J_{k}^{-1}J_{k} \Delta V_{k} - V^* \Vert_{\infty}\\
&= \Vert J_{k}^{-1}\left(J_k V_k - r(V_k) + r(V_k) + J_{k} \Delta V_{k} - J_k V^* \right)\Vert_{\infty}\\
&\leq \Vert J_{k}^{-1}\Vert_{\infty}\Vert J_k V_k - r(V_k) + r(V_k) + J_{k} \Delta V_{k} - J_k V^* \Vert_{\infty}\\
&\overset{(a)}{\leq} \frac{1}{1-\gamma} \left[ \Vert  r(V_k) - J_k \left( V_k - V^*  \right) \Vert_{\infty} \!+ \Vert J_k \Delta V_{k} + r(V_k) \Vert_{\infty}   \right]\\
&\overset{(b)}{\leq} \frac{1}{1-\gamma} \left[ \Vert  r(V_k) - J_k \left( V_k - V^*  \right) \Vert_{\infty} + \alpha_k\,\Vert r(V_k) \Vert_{\infty}   \right]\\
&= \frac{1}{1-\gamma} \left[ \Vert  r(V_k) - r(V^*) - J_k \left( V_k - V^*  \right) \Vert_{\infty}\right. \\
&\phantom{=} \left. \,\,\,\, +\, \alpha_k\,\Vert r(V_k) - r(V^*) \Vert_{\infty}   \right]\\
&\overset{(c)}{\leq} \frac{1}{1-\gamma} \left[ \Vert  r(V_k) - r(V^*) - J_k \left( V_k - V^*  \right) \Vert_{\infty}\right. \\
&\phantom{=} \left. \,\,\,\, +\, (1+\gamma)\alpha_k\,\Vert V_k - V^* \Vert_{\infty}   \right]\,,
\end{align*}
\endgroup
where $(a)$ follows from Lemma~\ref{lm:inverse_bound}, $(b)$ from the stopping condition~\eqref{eq: stopping_rule} and $(c)$ from Lemma~\ref{lm:lipschitz_const}.

Since $r$ is strongly semismooth at $V^*$, there exists a neighborhood $\mathcal{N}(V^*)$ such that, if $V_k \in \mathcal{N}(V^*)$, then
\begin{equation*}
\Vert V_{k+1} - V^* \Vert_{\infty} \leq \mathcal{O}\left( \Vert V_k - V^* \Vert_{\infty} ^2\right) + \frac{1+\gamma}{1-\gamma}\alpha_k \Vert V_k - V^*\Vert_{\infty}\,,
\end{equation*}
from which we can conclude local Q-linear convergence with rate $\frac{1+\gamma}{1-\gamma}\alpha<1$ if $\left\{\alpha_k \right\}\subseteq [0, \alpha]$ with $\alpha \in (0, \frac{1-\gamma}{1+\gamma})$ and local Q-superlinear convergence if $\lim_{k \rightarrow \infty}\alpha_k = 0$.~\qed  
\end{pf}

The results of Theorem~\ref{th:local_contraction} show that the local convergence properties of inexact policy iteration methods are strongly affected by the forcing sequence. In addition, because of the specific structure of the problem at hand, we can compute the range of $\alpha$-values for which local convergence is guaranteed as it solely depends on $\gamma$. 


\section{Numerical Evaluation}\label{sec:NE}
We evaluate the performance of iGMRES-PI against PI and OPI on an MDP with $10000$ states, $40$ actions and $\gamma=0.95$. The methods are implemented in Python using NumPy~\citep{harris2020array} and the simulations are run on an Intel(R) Core(TM) i7-10750H CPU @ 2.60GHz architecture. We enforce single-core execution for all algorithms to ensure a fair comparison.

\begin{figure}
\centering
\begin{subfigure}[b]{0.38\textwidth}
   \includegraphics[width=1\linewidth]{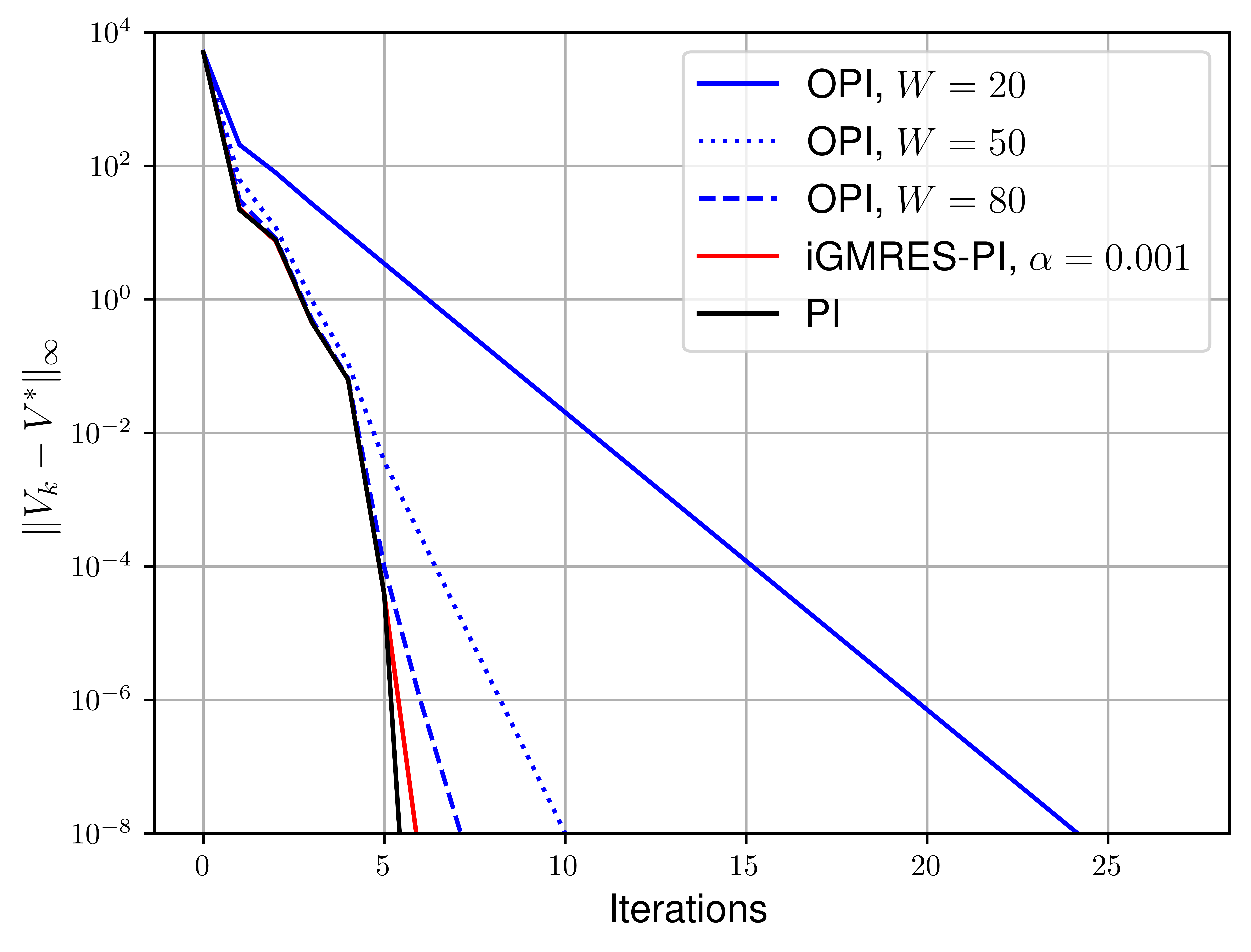}
\vspace*{-7mm}
   \caption{Infinity-norm of the suboptimality gap vs outer iterations for $\gamma=0.95$.}
   \label{fig:Ng1} 
\end{subfigure}

\begin{subfigure}[b]{0.38\textwidth}
   \includegraphics[width=1\linewidth]{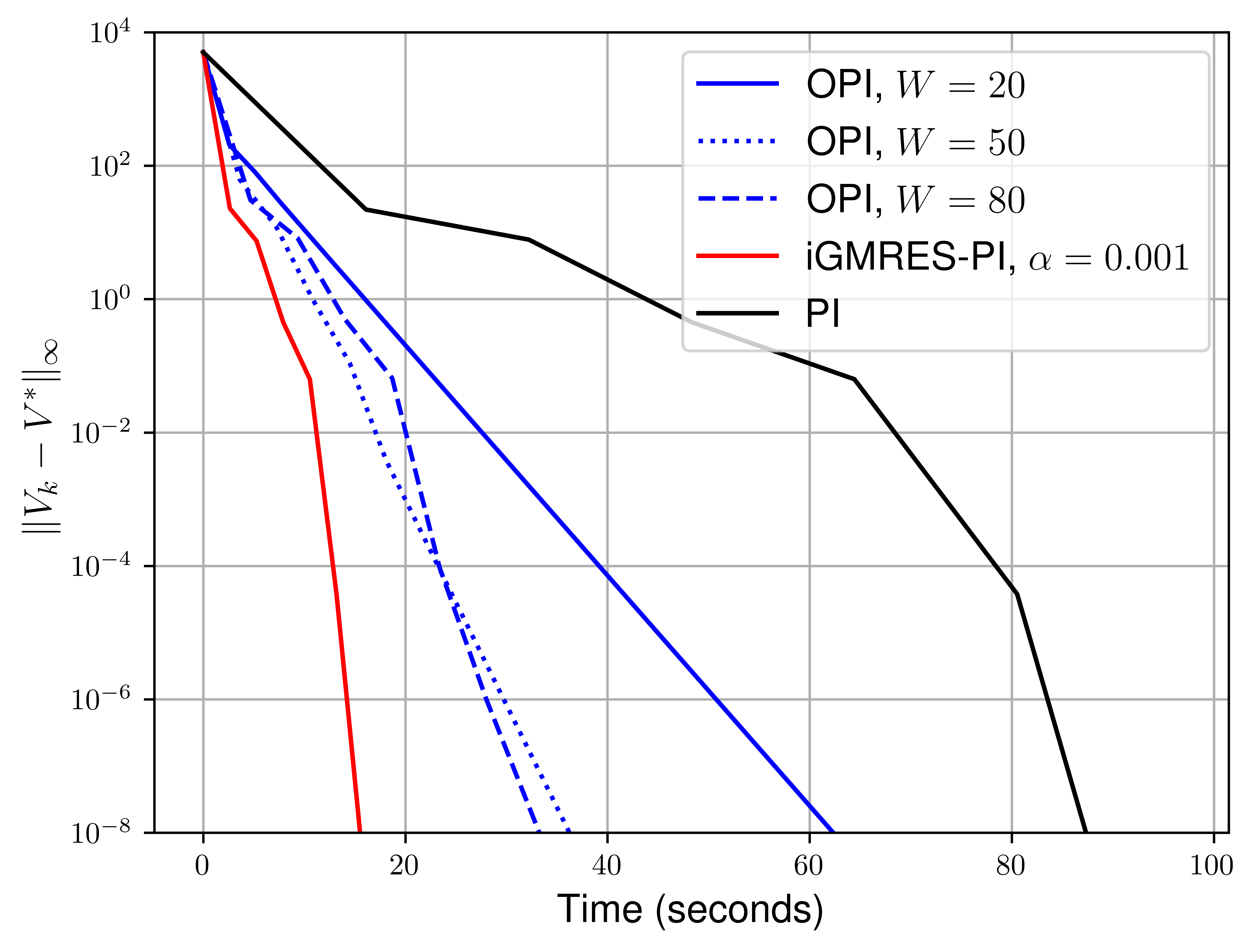}
   \vspace*{-7mm}
   \caption{Infinity-norm of the suboptimality gap vs time for $\gamma=0.95$.}
   \label{fig:Ng2}
\end{subfigure}
\begin{subfigure}[b]{0.38\textwidth}
   \includegraphics[width=1\linewidth]{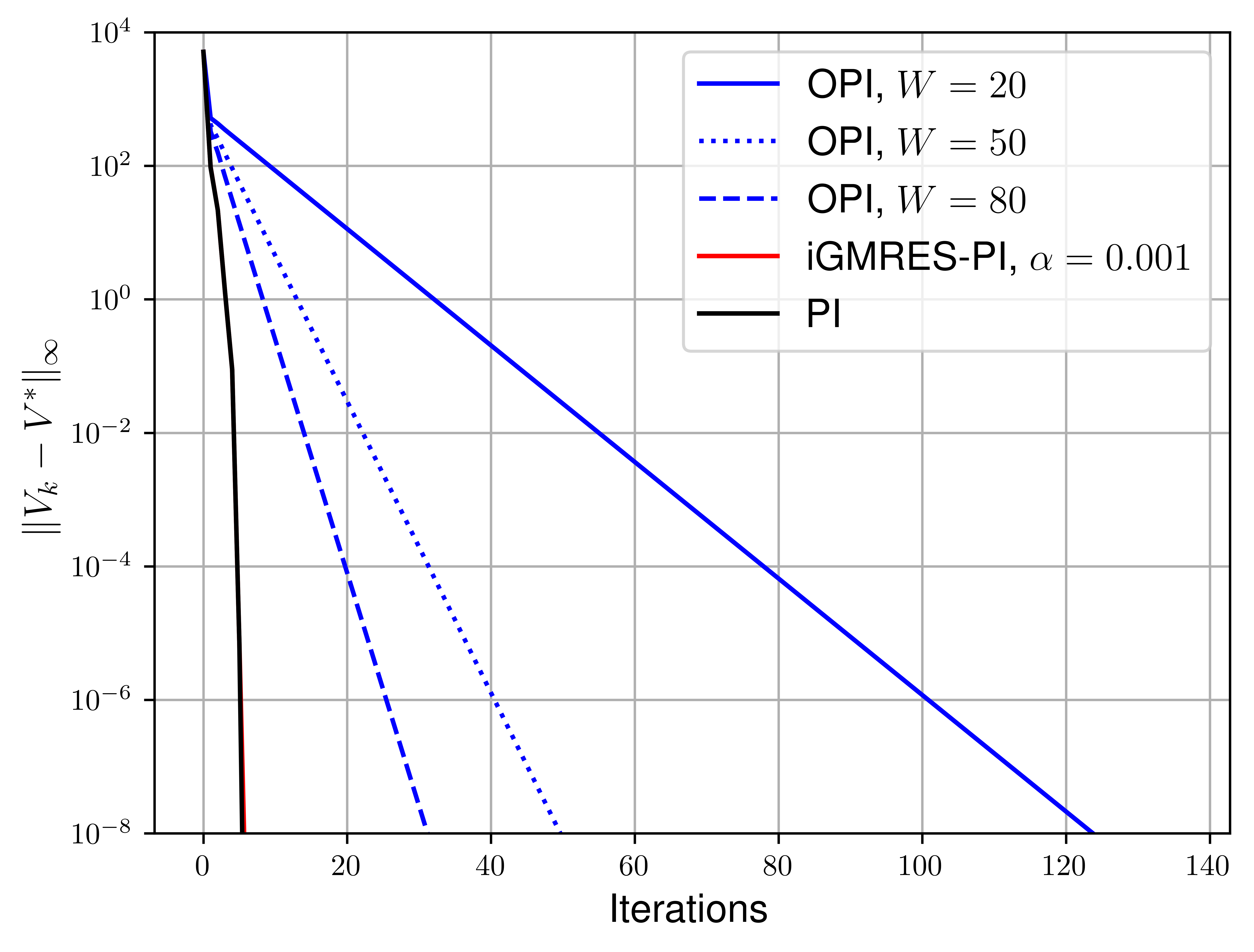}
   \vspace*{-7mm}   
   \caption{Infinity-norm of the suboptimality gap vs outer iterations for $\gamma=0.99$.}
   \label{fig:Ng3} 
\end{subfigure}

\begin{subfigure}[b]{0.38\textwidth}
   \includegraphics[width=1\linewidth]{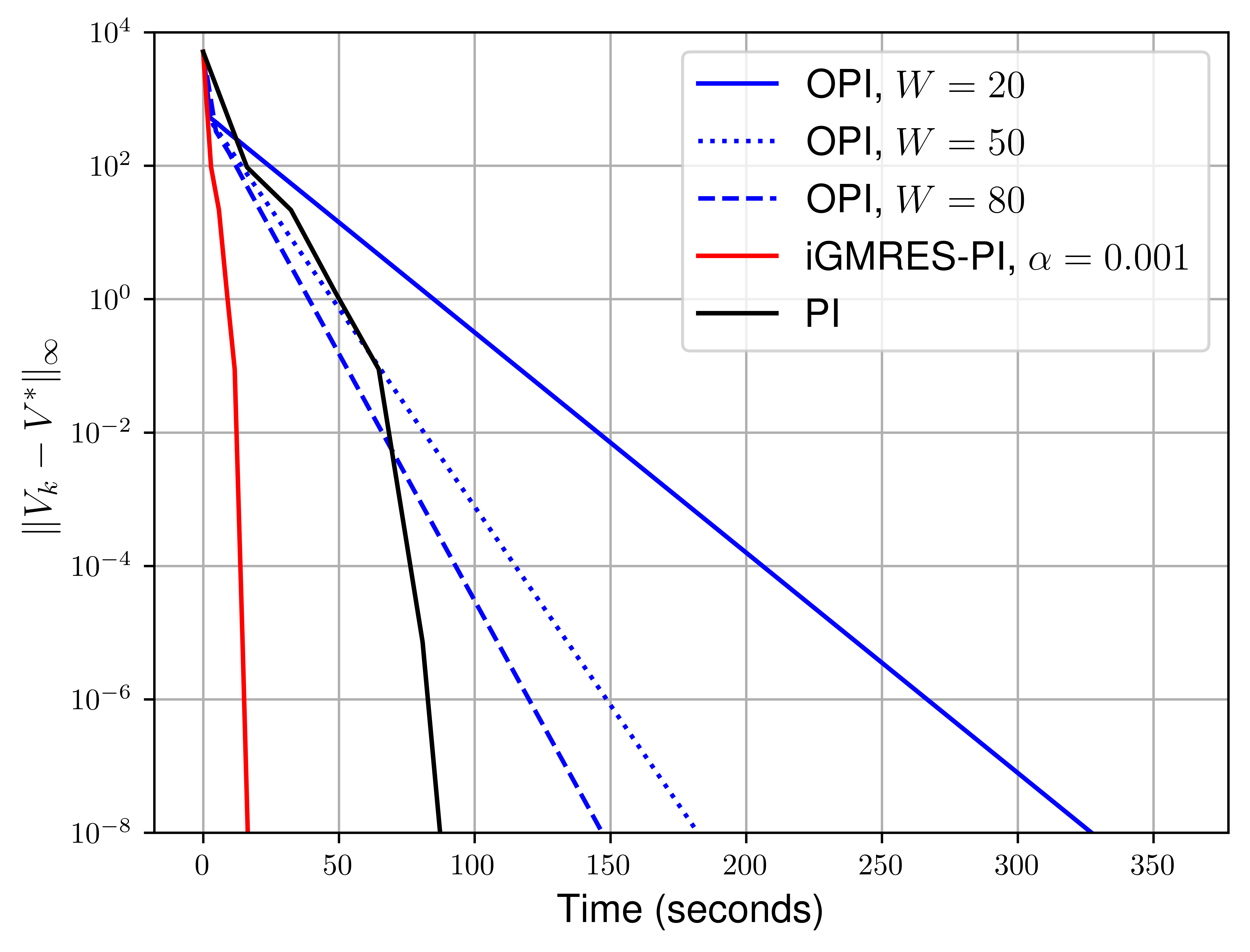}
   \vspace*{-7mm}
   \caption{Infinity-norm of the suboptimality gap vs time for $\gamma=0.99$.}
   \label{fig:Ng4}
\end{subfigure}
\vspace*{-3mm}
\caption{Performance of iGMRES-PI (red), PI (black) and OPI with a different number of inner iterations $W$ (blue) on an MDP with $n=10000$ and $m=40$.}
\label{fig:both}
\end{figure}

In Figure~\ref{fig:Ng1} we visualize the infinity-norm of the suboptimality gap versus the number of outer iterations.
As shown in Figure~\ref{fig:Ng1}, PI enjoys the fastest rate of convergence, followed by iGMRES-PI and OPI. As expected, the convergence rate of OPI improves by increasing the number of inner iterations $W$.
When considering time instead of outer iterations, the situation changes dramatically for PI. Its expensive iterations result in PI being the slowest converging method in terms of wall-clock time.  As shown in Figure~\ref{fig:Ng2}, PI takes $\sim 87$ seconds to reach convergence. A better trade-off between convergence rate and computational complexity is offered by OPI with $W=50$ and $W=80$, which achieves convergence in $\sim 36$ and $\sim 33$ seconds, respectively. Our iGMRES-PI greatly outperforms both PI and OPI, achieving convergence in only $\sim 15$ seconds and therefore attaining a $\times 5.8$ and $\times 2.2$ speedup with respect to PI and OPI, respectively.   

We then run the same benchmarks increasing the discount factor to 0.99. As depicted in Figures~\ref{fig:Ng3} and \ref{fig:Ng4}, this increase has a dramatic effect on the performance of OPI, while PI and iGMRES-PI's performance is essentially unaltered. These empirical results are in line with our observations in Figure~\ref{fig:GMRESvsVI_convergence}.

\section{Conclusions \& Future Work}

Taking inspiration from inexact semismooth Newton methods, we define a novel class of DP methods for large-scale applications which we call inexact policy iteration. We provide local contraction guarantees for the methods in this class and propose iGMRES-PI, an instance of inexact policy iteration based on the deployment of GMRES for the approximate policy evaluation step. 
We validate the performance superiority of iGMRES-PI against PI and OPI on a large-scale MDP.

Future work includes providing global convergence guarantees for inexact policy iteration, boosting the performance of GMRES via the design of an ad-hoc preconditioner for policy evaluation and the study of high-performance parallel and distributed variants of iGMRES-PI.

\bibliography{ifacconf}             
                                                   







\end{document}